%% file: newDifferentialsInHigherSpectralSequences10.tex
\documentclass[10pt]
{article}

\usepackage{amsthm}
\usepackage[margin=30mm]{geometry}
\usepackage{amsmath}%
\usepackage{amsfonts}%
\usepackage{amssymb}%
\usepackage{amscd}
\usepackage{graphicx}
\usepackage{color}
\usepackage{xcolor}
\usepackage{bbm}
\usepackage[all,cmtip]{xy}
\usepackage{rotating}
\usepackage{url}
\usepackage{enumitem}
\usepackage{subcaption}

\usepackage{eucal}
\DeclareMathAlphabet{\mathpzc}{OT1}{pzc}{m}{it}
\usepackage{mathrsfs}

\definecolor{grau}{rgb}{0.1,0.1,0.1}
\usepackage[colorlinks=true, urlbordercolor=grau, linkcolor=grau, citecolor=grau, urlcolor=grau]{hyperref}
\usepackage{breakurl}

\usepackage{tikz}
\tikzset{node distance=2cm, auto}

\usepackage{bibspacing}

\newtheorem{theorem}{Theorem}[section]

\newtheorem{lemma}[theorem]{Lemma}
\newtheorem{corollary}[theorem]{Corollary}

\newtheorem{observation}[theorem]{Observation}

\theoremstyle{definition}
\newtheorem{definition}[theorem]{Definition}
\newtheorem{example}[theorem]{Example}
\newtheorem{funfact}[theorem]{Fun fact}
\newtheorem{remark}[theorem]{Remark}
\newtheorem{remarks}[theorem]{Remarks}

\newcommand{\RR}{\mathbbm{R}} 
\newcommand{\QQ}{\mathbbm{Q}} 
\newcommand{\ZZ}{\mathbbm{Z}} 
\newcommand{\Zinf}{\overline{\mathbbm{Z}}} 
\DeclareMathOperator*{\wo}{\backslash}
\newcommand{\tof}[1]{\stackrel{#1}{\longrightarrow}} 
\newcommand{\Tof}[1]{\stackrel{#1}{\xrightarrow{\hspace*{1.7em}}}} 
\newcommand{\im}{\textnormal{im}} 
\DeclareMathOperator*{\coker}{coker}

\newcommand{\textdef}[1]{\textnormal{\textit{#1}}}
\newcommand{\iso}{\cong} 
\newcommand{\one}{\mathbbm{1}} 

\newcommand{\pt}{\textnormal{pt}} 
\newcommand{\incl}{\hookrightarrow} 
\newcommand{\st}{\ |\ } 
\newcommand{\Hom}{\textnormal{Hom}} 

\newcommand{\eps}{\varepsilon}

\newcommand{\id}{\textnormal{id}} 

\newcommand{\blank}{\underline{\ \ }}

\newcommand{\SL}{\textnormal{SL}}

\newcommand{\Tor}{\textnormal{Tor}}

\newcommand{\cat}{\mathcal}

\newcommand{\leqlex}{\leq_{\textnormal{lex}}} 
\newcommand{\leqlexsigma}[1]{\leq_{#1\textnormal{-lex}}}
\newcommand{\lesslexsigma}[1]{<_{#1\textnormal{-lex}}}

\newcommand{\Comp}{\textnormal{Comp}}
\newcommand{\concat}{*}
\newcommand{\tr}{t} 
\newcommand{\Sym}{\mathfrak{S}} 
\newcommand{\diffdir}{r} 
\newcommand{\diffdelta}{\delta} 
\newcommand{\filtind}{P} 
\newcommand{\upperbound}{u} 
\newcommand{\letters}{L} 
\newcommand{\transf}{M} 
\newcommand{\transfdelta}{T} 
\newcommand{\downset}{D} 
\newcommand{\lattice}{V} 
\newcommand{\letterext}{\textnormal{x}} 

\newcommand{\itref}[1]{\ref{#1}}

\begin{document}

\title{Higher spectral sequences}

\author{%
Benjamin Matschke%
\setcounter{footnote}{-1}%
\\
\small Boston University\\
\small matschke@bu.edu
}
\date{July 6, 2021}

\maketitle

\begin{abstract}
In this article we construct what we call a higher spectral sequence for any chain complex (or topological space)
that is filtered in $n$ compatible ways.
For this we extend the previous spectral system construction of the author, and we show that it admits considerably more differentials than what was previously known.
%
As a result, this endows the successive Leray--Serre, Grothendieck, chromatic--Adams--Novikov, and Eilenberg--Moore spectral sequences of the author with the structure of a higher spectral sequence.
Another application is a universal coefficient theorem analog for spectral sequences.
\end{abstract}

%
%
%
%

\renewcommand\contentsname{Contents}
\setcounter{tocdepth}{1}

\section{Introduction}

Consider either 
\begin{enumerate}[label={(\alph*)}]
\item\label{itFilteredC}
a chain complex $(C,d)$ that is filtered in $n$ different ways over the integers, or

\item\label{itFilteredX}
a topological space (or a spectrum) $X$ that is filtered in $n$ different ways over the integers, together with a generalized homology theory $h_*$.
\end{enumerate}
In the classical case $n=1$ we obtain a spectral sequence with limit\footnote{in the sense of~\cite{Mat13succSpectralSequences}} $H_*(C)$ and $h_*(X)$, respectively.
In the general case $n\geq 1$, 
under a natural compatibility assumption of the $n$ filtrations,
we will construct what we call a \emph{higher spectral sequence}.
For its construction, we extend the spectral system construction of~\cite{Mat13succSpectralSequences}. 
%
While it has the same limit $H_*(C)$ and $h_*(X)$, respectively, it contains considerably more pages than the $n$ naturally associated spectral sequences.
The pages considered in this paper are indexed over admissible words, and they can be related via differentials, group extensions, and natural isomorphisms. 
\noindent
More precisely, given $C$ or $X$ as in \itref{itFilteredC} or \itref{itFilteredX} above, we construct for any admissible word $\omega\in \letters^*_a$ (see Definition~\ref{defAdmissibleWord}) over the alphabet
\(
\letters = \{1,\ldots,n,1^\infty,\ldots,n^\infty,\letterext\}
\)
a so-called $\omega$-page, which is a collection of Abelian groups 
\[
S(P;\omega).
\]
Here $P$ ranges over a quotient $\ZZ^n/V_\omega\iso\ZZ^{n-k}$, where $k$ is the number of letters $\letterext$ in~$\omega$.
In the alphabet~$\letters$, a letter $j\in[n]$ stands for taking homology with respect to the $j$'th differential, $j^\infty$ denotes the same but infinitely often, and $\letterext$ stands for a group extension process.
The word $\omega$ is read from left to right and it tells how $S(P;\omega)$ is ``constructed'' from the first page, which itself is indexed by the empty word~$\eps$.

In ordinary spectral sequences we have $n=1$, $L=\{1,1^\infty,\letterext\}$, and for $\omega=1^{r-1}$ ($r\geq 1$) the $\omega$-page consists of the columns in~$E^r_{**}$, which are indexed over $P\in\ZZ$. 
The letter $1$ stands for the connection between any $E^r_{**}$ and~$E^{r+1}_{**}$, $1^\infty$ stands for the connection between any $E^r_{**}$ and $E^\infty_{**}$, and $\letterext$ for the connection between $E^\infty_{**}$ and the ``limit'' of the spectral sequence, e.g.\ $H_*(C)$ and $h_*(X)$ in the above settings.


The analog of this for higher spectral sequences can be informally summarized as follows.
\begin{enumerate}
\item \label{itIntroFirstPage} 
The first page $S(P;\eps)$, $P\in\ZZ^n$, consists of the homologies of the smallest pieces (subquotients of $C$, and pairs of subsets of $X$) that can be constructed from the common refinement of the $n$ given filtrations. 

\item \label{itIntroDifferentials} 
At most of the pages $S(P;\omega)$, there exist $n$ differentials (unless they are saturated), which are indexed by $[n] = \{1,\ldots,n\}$.
Taking homology with respect to the $j$'th differential yields $S(P;\omega * j)$.
Here, ``*'' denotes concatenation of words.

\item \label{itIntroInfinitelyManyDifferentials} 
One can take homology in one direction $j$ an infinite number of times to arrive at $S(P;\omega*j^\infty)$. This saturates the differential with index~$j$. 

\item \label{itIntroExtensions} 
After any such step, we perform a group extension process to arrive at $S(P;\omega*j^\infty\letterext)$.

\item \label{itIntroLimit} 
If $j^\infty\letterext$ appears in $\omega$ for each $j\in [n]$, we call $\omega$ final and $S(P;\omega)$ equals $H_*(C)$ or $h_*(X)$, respectively.
\end{enumerate}

This can be used as follows.
Step~\itref{itIntroDifferentials} gives a connection between the first page (Step~\itref{itIntroFirstPage}) and arbitrary $\omega$-pages for $\omega\in[n]^*$.
%
%
We can proceed with Step~\itref{itIntroInfinitelyManyDifferentials} and take homology infinitely often in one direction $j\in [n]$.
Note that as with usual spectral sequences, $S(\filtind;\omega\concat j^\infty)$ may indeed be a proper subquotient of the limit of $S(P;\omega\concat j^k)$ as $k\to\infty$, compare with Weibel~\cite{Wei94homologicalAlgebra}, Boardman~\cite{Boardman99conditionallyConvergentSS}, McCleary~\cite{McC01userGuideToSS}.
Then we proceed with Step~\itref{itIntroExtensions}, which connects to $S(\filtind;\omega\concat j^\infty \letterext)$.
As with the extension process in usual spectral sequences, the filtration $(F_i)$ of $S(\filtind;\omega\concat j^\infty \letterext)$ in Theorem~\ref{thmDifferentials}\itref{itMainThmExtensions} may be neither Hausdorff nor exhaustive, and even if it is, $S(\filtind;\omega\concat j^\infty \letterext)$ may not be complete with respect to $(F_i)$.
As usual these two problems in Steps~\itref{itIntroInfinitelyManyDifferentials} and~\itref{itIntroExtensions} can be serious obstacles for computations, but they are the standard ones in ordinary spectral sequences.
Arriving at $S(\filtind;\omega\concat j^\infty \letterext)$ we can start again at Step~\itref{itIntroDifferentials} and iterate until $\omega$ is final.
At that point, Step~\itref{itIntroLimit}, we arrived at the limit of the higher spectral sequence.


\paragraph{Outline.}
This paper is organized as follows.
In Section~\ref{secPreliminaries} we review the notions and basic properties of exact couple systems over $D(\ZZ)$ (``higher exact couples'') and their associated spectral systems (``higher spectral sequences'').
This is a common formal framework that includes both major special cases \itref{itFilteredC} and \itref{itFilteredX} from above.

In Section~\ref{secOmegaPages} the $\omega$-pages are defined.
In Section~\ref{secMainTheorem} the main theorem is stated, which makes the above informal description of higher spectral sequences and their basic properties precise.
Some examples 
are given: Higher Leray--Serre, Eilenberg--Moore, and Grothendieck spectral sequences, as well as two universal coefficient higher spectral sequences.
%
The main theorem is proved in Section~\ref{secProofOfMainTheorem}.
In Section~\ref{secCaseN2} we discuss some more properties of the $n=2$ case. 

\paragraph{Acknowledgements.}

This work was supported by NSF Grant DMS-0635607 at Institute for Advanced Study, by an EPDI fellowship at Institut des Hautes \'Etudes Scientifiques,  Forschungsinstitut f\"ur Mathematik (ETH Z\"urich), and the Isaac Newton Institute for Mathematical Sciences, by Max-Planck-Institute for Mathematics Bonn, and by Simons Foundation grant {\#}550023 at Boston University (in chronological order).%

\section{Preliminaries on spectral systems} \label{secPreliminaries}

We recall the necessary background from~\cite{Mat13succSpectralSequences}.
Let $n\geq 1$ and $[n]:=\{1,\ldots,n\}$. 
Let $e_1,\ldots,e_n$ be the standard basis vectors in $\ZZ^n$, and $\one:=(1,\ldots,1)^\tr\in\ZZ^n$.
$\ZZ^n$ is a poset via $(x_1,\ldots,x_n)\leq (x'_1,\ldots,x'_n)$ if and only if $x_i\leq x_i'$ for all~$i$.

Throughout the paper, let $I:=D(\ZZ^n)$ denote the lattice of downsets of~$\ZZ^n$. (Everything in this paper can also be done for filtrations over $D(\Zinf^n)$, where $\Zinf=\ZZ\cup\{\pm\infty\}$; here we consider only $D(\ZZ^n)$ because it makes the presentation cleaner.)

$I$ has minimum $-\infty:=\emptyset$ and maximum $\infty:=\ZZ^n$.
We write $I_k:=\{(p_1,\ldots,p_k)\in I^k\st p_1\geq\ldots\geq p_k\}$, which is again a poset via $(p_1,\ldots,p_k)\leq (p'_1,\ldots,p'_k)$ if and only if $p_i\leq p_i'$ for all~$i$.

For us a chain complex is an Abelian group $C$ together with an endomorphism $d:C\to C$ with $d\circ d=0$, and its homology is $H(C,d):=\ker(d)/\im(d)$; the grading is not of importance for us.
An $I$-filtration of $C$ is a family of subchain complexes $(F_p)_{p\in I}$ such that $F_q\subseteq F_p$ whenever $q\leq p$.

Similarly, if $X$ is a topological space then an $I$-filtration of $X$ is a family of (usually open) subspaces $(X_p)_{p\in I}$ such that $X_q\subseteq X_p$ whenever $q\leq p$.

Whenever we have an $I$-filtered chain complex $(C,d)$, or an $I$-filtered space $X$ together with a generalized homology theory $h_*$, we can associate a so-called exact couple system via
\begin{equation}
\label{eqExactCoupleSystemOfFilteredChainComplex}
E^p_q:=H(F_p/F_q)
\end{equation}
or
\begin{equation}
\label{eqExactCoupleSystemOfFilteredSpace}
E^p_q:=h_*(X_p,X_q),
\end{equation}
respectively, which is defined as follows.

\begin{definition}[Exact couple system]
And \textdef{exact couple system over $I$} is a collection of Abelian groups $(E^p_q)_{(p,q)\in I_2}$ together with homomorphisms $\ell^{p,q}_{p',q'}:E^p_q\to E^{p'}_{q'}$ for any $(p,q)\leq (p',q')$ and homomorphisms $k_{p,q}:E^p_q\to E^q_{-\infty}$ for any $(p,q)\in I_2$, such that the following properties are satisfied:
\begin{enumerate}
\item $\ell^{p',q'}_{p'',q''}\circ\ell^{p,q}_{p',q'}=\ell^{p,q}_{p'',q''}$.
\item The triangles
\[
\xymatrix{
E^q_{-\infty} \ar[rr]^{\ell^{q,-\infty}_{p,-\infty}} &                        & E^p_{-\infty} \ar[dl]^{\ell^{p,-\infty}_{p,q}} \\
                     & E^p_q \ar[ul]^{k_{pq}} &
}
\]
are exact.
\item The diagrams
\[
\xymatrix{
E^p_q \ar[d]_{\ell^{pq}_{p'q'}} \ar[r]^{k_{pq}} & E^q_{-\infty} \ar[d]^{\ell^{q,-\infty}_{q',-\infty}} \\
E^{p'}_{q'} \ar[r]_{k_{p'q'}}                   & E^{q'}_{-\infty}
}
\]
commute.
\end{enumerate}
\end{definition}

Let $E$ be an exact couple system over $I$.
There is a natural differential $d_{pqz}:E^p_q\to E^q_z$ for any $(p,q,z)\in I_3$ defined by $d_{pqz}:=\ell^{q,-\infty}_{q,z}\circ k_{pq}$.
With this we define an associated \emph{spectral system over $I$} via
\begin{equation}
\label{eqStermsOfExactCoupleAnalog}
S^{pz}_{bq}:=\frac{\ker(d_{pqz}:E^p_q\to E^q_z)}{\im(d_{bpq}:E^b_p\to E^p_q)},\quad (b,p,q,z)\in I_4.
\end{equation}
In the special case $I=D(\ZZ^n)$, we call the spectral system $S$ the \emph{higher spectral sequence of $E$}.
At a first glance this is just a collection of Abelian groups, one for each element in $I_4$, however there are many connections between them:

First note that the usual goal of computation, $E^\infty_{-\infty}$, appears as $S^{\infty,-\infty}_{\infty,-\infty}$.
It is called the \emph{limit} of this spectral system (this is just a name; it does not imply any convergence or comparison property).
Moreover, terms of the form $S^{pq}_{pq}=E^p_q$ are usually known when $p$ covers $q$, that is, when $|p\wo q|=1$.
The collection of these terms is called the \emph{first page} of~$S$.

The following facts are proved in~\cite{Mat13succSpectralSequences}.
For any $(b,p,q,z)\leq (b',p',q',z')$ in $I_4$, $\ell^{pq}_{p'q'}$ induces maps
\[
S^{pz}_{bq}\to S^{p'z'}_{b'q'},
\]
which we call maps induced by inclusion.
When there is no confusion, we abbreviate all of them as~$\ell$.

\begin{lemma}[Extensions]
\label{lemExtensionProperty}
For any $z\leq p_1\leq p_2\leq p_3\leq b$ in $I$, we have a short exact sequence of maps induced by inclusion,
\begin{equation}
\label{eqExtensionProperty}
0\to S^{p_2,z}_{b,p_1}\to S^{p_3,z}_{b,p_1}\to S^{p_3,z}_{b,p_2}\to 0.
\end{equation}
\end{lemma}

\begin{lemma}[Differentials]
\label{lemDifferentials}
For any $(b,p,q,z),(b',p',q',z')\in I_4$ with $z\leq p'$ and $q\leq b'$ there are natural differentials
\begin{equation}
\label{eqDifferential}
d: S^{pz}_{bq}\to S^{p'z'}_{b'q'},
\end{equation}
which commute with $\ell$, that is, $\ell\circ d = d\circ \ell$.
\end{lemma}

\begin{lemma}[Kernels and cokernels]
\label{lemKernelsAndCokernels}
For any $(b,p,q,z),(b',p',q',z')\in I_4$ with $z=p'$ and $q=b'$ we have
\[
\ker\left(d: S^{pz}_{bq}\to S^{p'z'}_{b'q'}\right) = S^{pq'}_{bq}
\]
and
\[
\coker\left(d: S^{pz}_{bq}\to S^{p'z'}_{b'q'}\right) = S^{p'z'}_{pq'}
\]
\end{lemma}

\begin{lemma}[$\infty$-page as filtration quotients]
\label{lemExactCouplesInfinityPageFiltrationQuotients}
$E^\infty_{-\infty}$ can be $I$-filtered by
\[
G_p:=\im (\ell:E^p_{-\infty}\to E^\infty_{-\infty}) \iso S^{p,-\infty}_{\infty,-\infty},\ \ \ p\in I.
\]
Furthermore the $S$-terms on the $\infty$-page are filtration quotients
\[
S^{p,-\infty}_{\infty,q}\iso G_p/G_q.
\]
\end{lemma}

\begin{lemma}[$\infty$-page as quotient kernels]
\label{lemExactCouplesInfinityPageQuotientKernels}
$E^\infty_{-\infty}$ has quotients
\[
Q_p:=\frac{E^\infty_{-\infty}}{\ker(\ell:E^\infty_{-\infty}\to E^\infty_p)}\iso S^{\infty,-\infty}_{\infty,p},\ \ \ p\in I.
\]
Furthermore the $S$-terms on the $\infty$-page are quotient kernels
\[
S^{p,-\infty}_{\infty,q}\iso \ker(Q_q\to Q_p).
\]
\end{lemma}

\begin{definition}[Excision]
An exact couple system $E$ over $I$ is called \textdef{excisive} if for all $a,b\in I$,
\[
E^a_{a\cap b} \tof{\ell} E^{a\cup b}_b
\]
is an isomorphism.
\end{definition}

The exact couple system~\eqref{eqExactCoupleSystemOfFilteredSpace} is automatically excisive by the excision axiom of $h_*$ if the subspaces $X_p$ are all open.
Note however that~\eqref{eqExactCoupleSystemOfFilteredChainComplex} is in general not excisive, though in many applications it is, for example when $C=\bigoplus_{P\in\ZZ^n}C_p$ (as Abelian group) and $(F_p)_I$ is the canonical $I$-filtration given by $F_p=\bigoplus_{P\in p}C_p$.
More generally, \eqref{eqExactCoupleSystemOfFilteredChainComplex} is excisive if the filtration $(F_p)_I$ is distributive in the sense that $F_{a\cap b} = F_a\cap F_b$ and $F_{a\cup b} = F_a + F_b$ for all $a,b\in I$.

In the settings \itref{itFilteredC} and \itref{itFilteredX} from the introduction,
we refine the $n$ given filtrations to obtain $I$-filrations of $C$ and $X$, respectively, by taking unions of $n$-wise intersections of the filtration pieces (see\
\cite[\S{3}]{Mat13succSpectralSequences}).
This induces associated exact couple systems~\eqref{eqExactCoupleSystemOfFilteredChainComplex} and~\eqref{eqExactCoupleSystemOfFilteredSpace}, respectively.
The natural compatibility assumption on the $n$ filtrations as mentioned in the introduction refers precisely to the excision property of these exact couple systems.

\medskip

Let us think of $J:=\ZZ^n$ as an undirected graph, whose vertices are the elements of $J$, and $x,y\in J$ are adjacent if they are related, i.e. $x\geq y$ or $x\leq y$ (coordinate-wise).
For $(b,p,q,z)\in I_4$, let $Z(z,q,p,b)\subseteq J$ denote the union of all connected components of $p\wo z$ that intersect $p\wo q$, and let $B(z,q,p,b)\subseteq I$ denote the union of all connected components of $b\wo q$ that intersect $p\wo q$.

\begin{lemma}[Natural isomorphisms]
\label{lemNaturalIsomorphisms2}
In an excisive exact couple system $E$ over $I=D(J)$, $S^{pz}_{bq}$ is uniquely determined up to natural isomorphism by $Z:=Z(z,q,p,b)$ and $B:=B(z,q,p,b)$.
\end{lemma}

We also write $S^Z_B$ for $S^{pz}_{bq}$, which is only defined up to natural isomorphisms.
A word of warning:
This $B$-$Z$-description of $S^{pz}_{bq}$ looks quite appealing.
However it may be combinatorially non-trivial to check whether some given $B$ and $Z$ come from some $(b,p,q,z)$, and if so there might be several good choices.
Moreover, it can be quite challenging to see whether there is a differential from $S^{Z_1}_{B_1}$ to $S^{Z_2}_{B_2}$ and what the resulting kernels and cokernels are in this case.


\section{Higher spectral sequences} \label{secDifferentials}

Throughout this section let us fix an excisive exact couple system $E$ over $I=D(\ZZ^n)$.

\subsection{Definition of $\omega$-pages}
\label{secOmegaPages}

Define an alphabet $\letters$,
\[
\letters:=\{1,\ldots,n,1^\infty,\ldots,n^\infty,\letterext\}.
\]

\begin{remark}[Some intuition]
Here, a letter $j\in[n]$ stands for taking homology with respect to the $j$'th differential, $j^\infty$ denotes the same but infinitely often, and $\letterext$ stands for a group extension process.
In ordinary spectral sequences, $n=1$, and the letter $1$ stands the connection between some $E^r_{**}$ and $E^{r+1}_{**}$, $1^\infty$ stands for the connection between some $E^r_{**}$ and $E^\infty_{**}$, and $\letterext$ for the connection between $E^\infty_{**}$ and the ``limit'' of the spectral sequence, e.g. $H(C)$ if the spectral sequence comes from a $\ZZ$-filtration of a chain complex $C$.
\end{remark}

Let $\letters^*$ denote the monoid of words of finite length with letters in $\letters$.
Denote the empty word by $\eps$, the concatenation of two words $\omega$ and $\omega'$ by $\omega\concat\omega'$, $\omega^n:=\omega\concat\ldots\concat\omega$ ($n$ times), and the length of $\omega$ by $|\omega|$.
$L^*_a$ becomes a poset via $\tau\leq\omega$ if and only if $\tau$ is a prefix of $\omega$, that is, a subword that starts from the beginning ($\tau=\eps$ and $\tau=\omega$ are allowed).

\begin{definition}[Admissible words]
\label{defAdmissibleWord}
Call a finite word $\omega\in \letters^*$ \textdef{admissible} if the following holds:
\begin{enumerate}
\item if $j^\infty$ appears, the subsequent subword of $\omega$ contains neither $j$ nor $j^\infty$,
\item the only letter allowed directly after $j^\infty$ is $\letterext$,
\item any $\letterext$ occurring in $\omega$ comes directly after some $j^\infty$.
\end{enumerate}
If furthermore $\omega$ contains subwords $j^\infty \letterext$ for all $j\in[n]$ then $\omega$ is called \textdef{final}.
\end{definition}

An exemplary final word for $n=3$ is $123122^\infty\letterext 133313^\infty\letterext 111^\infty\letterext$ and any prefix of a final word is admissible.
Let $\letters^*_a$ denote the set of all admissible words in $\letters^*$.
Define $X(\omega)\subseteq[n]$ as the set of $j\in[n]$ such that $j^\infty\letterext$ is a subword of~$\omega$, and $Y(\omega):=[n]\wo X(\omega)$.
$X(\omega)$ is so to speak the set of saturated indices along which the extension process has been already made, and $Y(\omega)$ is the set of unsaturated indices along which we still have differentials.

For $\omega\in \letters^*_a$, $i,j\in[n]$, we inductively define $\diffdir_\omega^i,\diffdelta_\omega^i\in\ZZ^n$ and $B_\omega, Z_\omega \subset\ZZ^n$ as follows.
Put $\diffdir_\eps^i:=e_i$, $\diffdir^i_{\omega\concat j^\infty}:=\diffdir^i_\omega$, $\diffdir^i_{\omega\concat \letterext}:=\diffdir^i_\omega$, and
\[
\diffdir_{\omega\concat j}^i:=
\begin{cases}
\diffdir_\omega^i                     & \mbox{if } i\neq j, \\
\diffdir_\omega^i+\diffdelta_\omega^i & \mbox{if } i=j,
\end{cases}
\]
where $\diffdelta_\eps^i:=e_i$, $\diffdelta^i_{\omega\concat j^\infty}:=\diffdelta^i_\omega$, $\diffdelta^i_{\omega\concat\letterext}:=\diffdelta^i_\omega$, and
\[
\diffdelta_{\omega\concat j}^i:=
\begin{cases}
\diffdelta_\omega^i                     & \mbox{if } i\in X(\omega)\cup\{j\}, \\
\diffdelta_\omega^i-\diffdelta_\omega^j & \mbox{if } i\in Y(\omega)\,\wo\,\{j\}.
\end{cases}
\]
For $\omega\in[n]^*$, $\diffdelta_\omega^i=\one-\sum_{k\in[n]\wo i}\diffdir_\omega^k$.
\begin{remark}[Some intuition 2]
$\diffdir^i_\omega$ will be the negated direction of the $i$'th differential at the $\omega$-page,
and $\diffdelta^i_\omega$ is the negated change of direction for the $i$'th differential that occurs when taking homology with respect to it.
In ordinary spectral sequences, $n=1$, and for $\omega=1^{r-1}$ the $\omega$-page consists of the columns in $E^r_{**}$, with $\diffdir^1_\omega=r$, and $\diffdelta^1_\omega = 1$.
\end{remark}

\begin{figure}[htb]
\centerline{
\input{differentialTree2D3.pspdftex}
} 
\caption{All $B_\omega$ with $|\omega|\leq 4$, $\omega_1=1$, and $n=2$.
For each $B_\omega$, the origin is marked with a solid square, and the two points $\diffdir^i_\omega-e_i/2$ are marked with a black dot.}
\label{figB2D}
\end{figure}
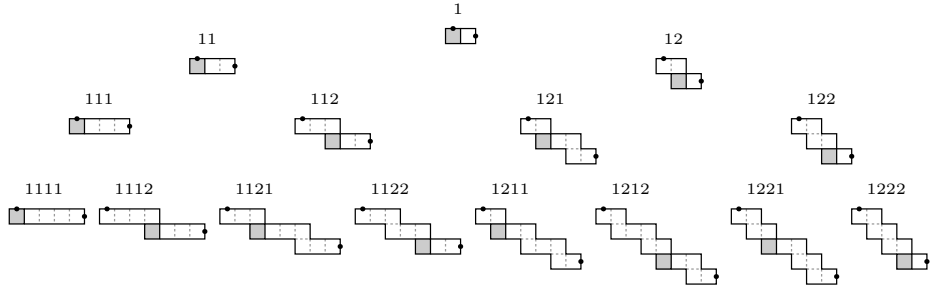

\begin{figure}[htb]
  \centering
  \begin{minipage}[b]{0.15\textwidth}
	\centering
	\includegraphics[scale=0.25]{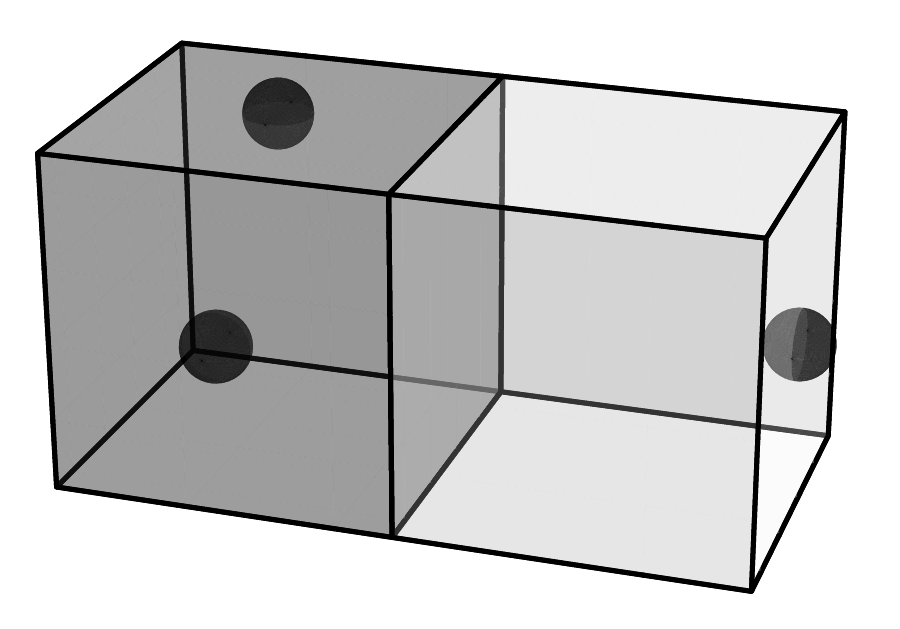}
    \subcaption*{$B_{1}$}
  \end{minipage}
\hfill
  \begin{minipage}[b]{0.15\textwidth}
	\centering
	\includegraphics[scale=0.25]{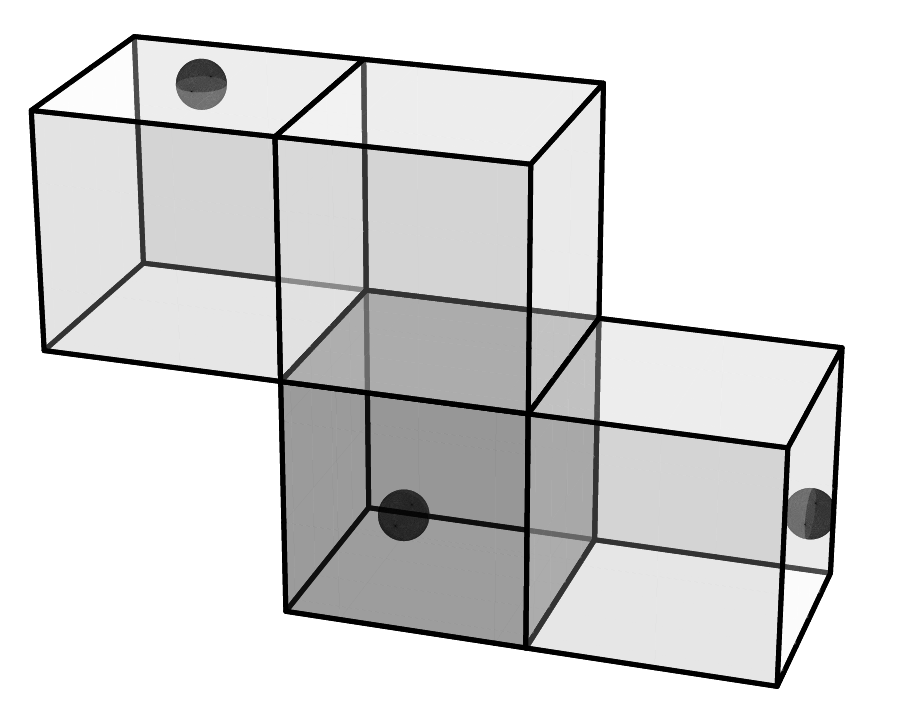}
    \subcaption*{$B_{12}$}
  \end{minipage}
\hfill
  \begin{minipage}[b]{0.15\textwidth}
	\centering
	\includegraphics[scale=0.25]{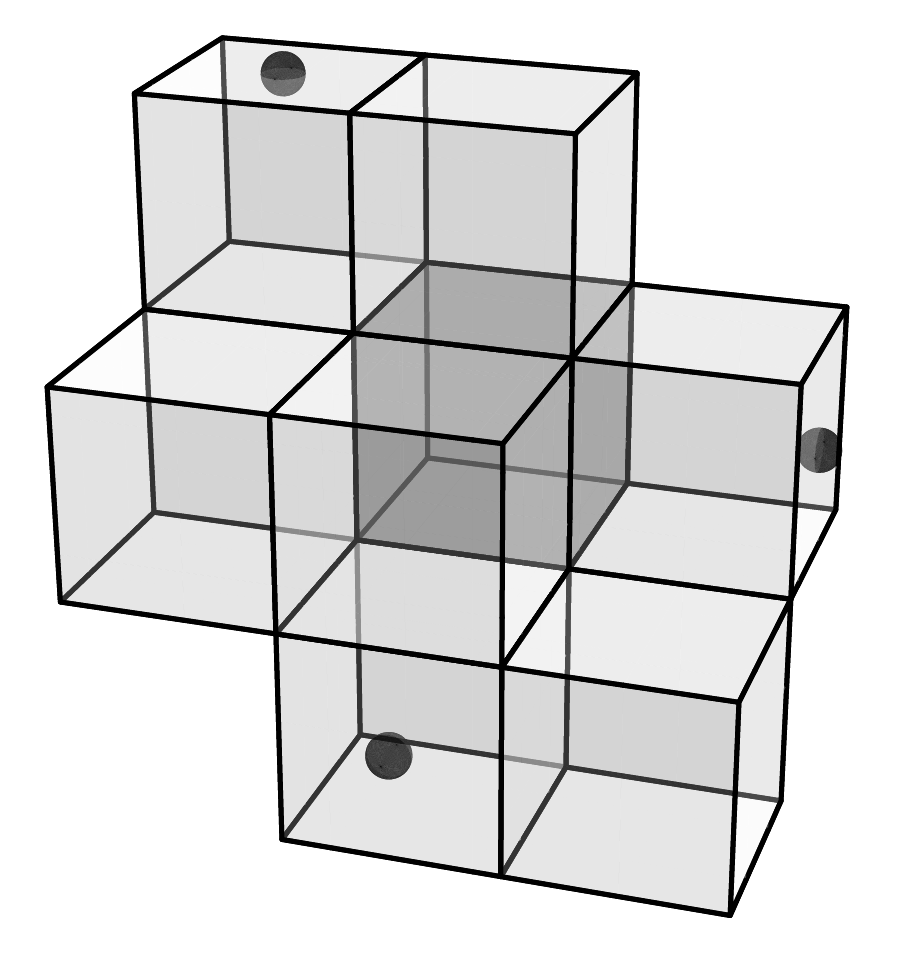}
    \subcaption*{$B_{123}$}
  \end{minipage}
\hfill
  \begin{minipage}[b]{0.15\textwidth}
	\centering
	\includegraphics[scale=0.25]{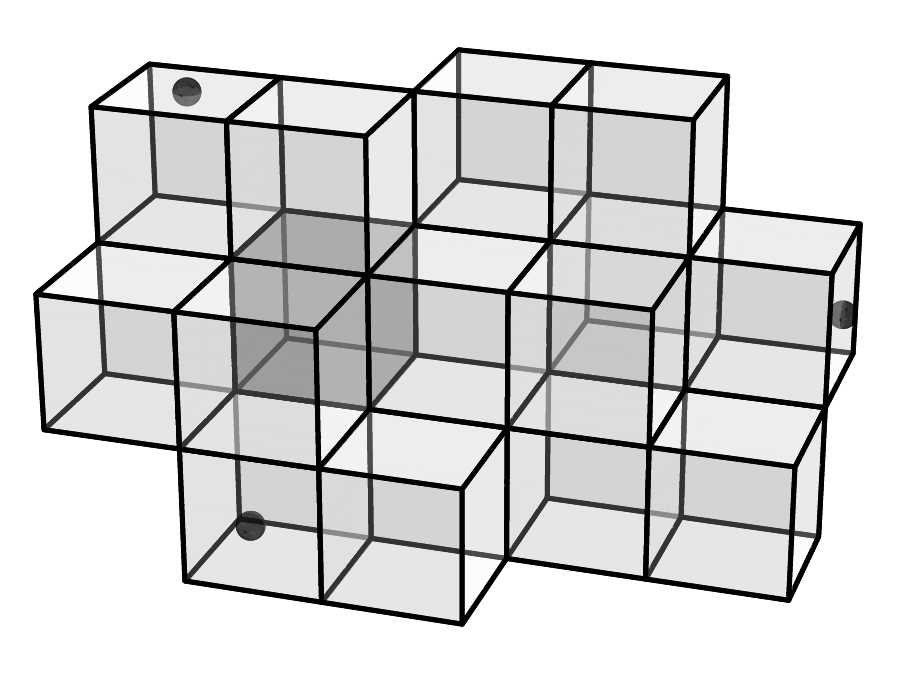}
    \subcaption*{$B_{1231}$}
  \end{minipage}
\hfill
  \begin{minipage}[b]{0.15\textwidth}
	\centering
	\includegraphics[scale=0.25]{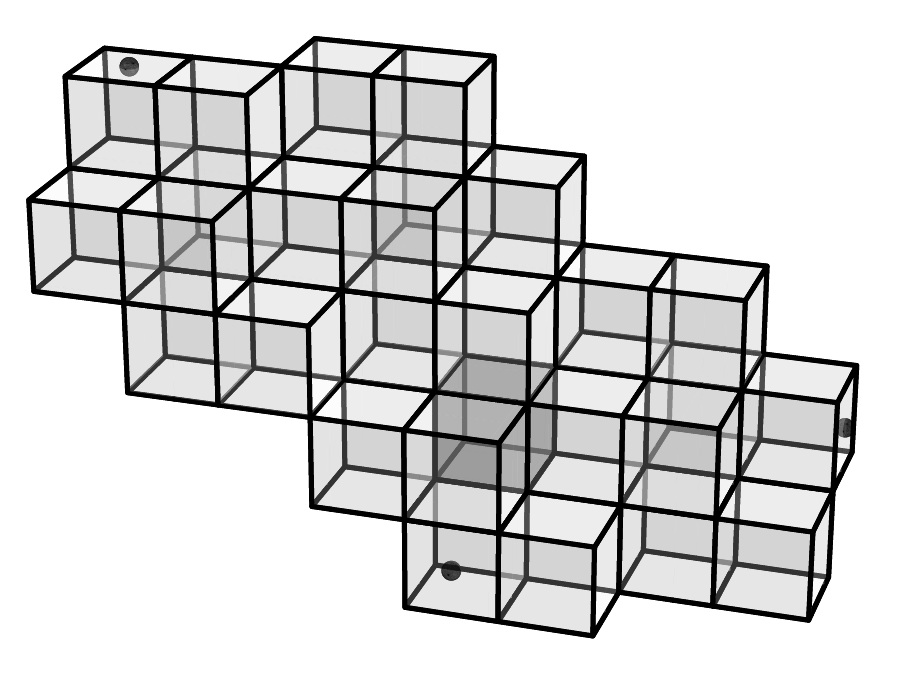}
    \subcaption*{$B_{12312}$}
  \end{minipage}
\hfill
  \begin{minipage}[b]{0.15\textwidth}
	\centering
	\includegraphics[scale=0.25]{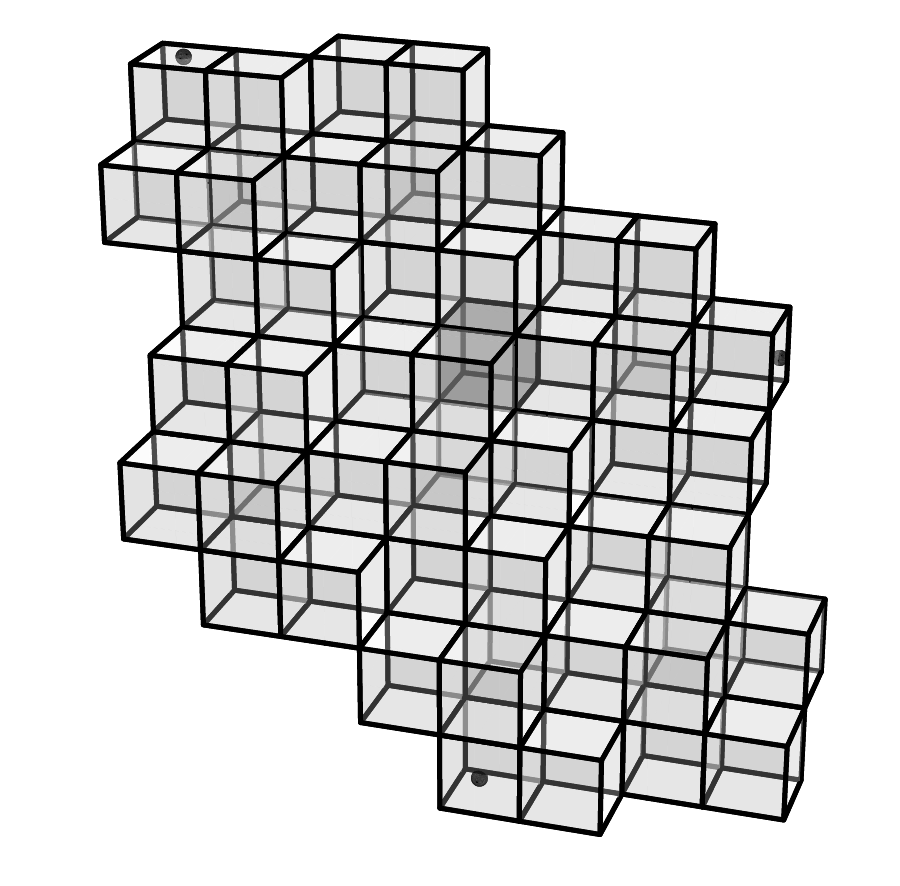}
    \subcaption*{$B_{123123}$}
  \end{minipage}
\caption{$B_\omega$ for $\omega=1, \ldots, 123123$ and $n=3$. For each $B_\omega$, the origin is marked with a dark cube, and the three points $\diffdir^i_\omega-e_i/2$ with a black dot.}
\label{figB3D}
\end{figure}

Further put $B_\eps:=\{0\}$,
\begin{align}
\label{eqBasZonotope}
B_{\omega\concat j}:=\ & B_\omega + \{0,\diffdelta^j_\omega\}, \\
B_{\omega\concat j^\infty}:=\ & B_\omega + \ZZ_{\geq 0}\cdot \diffdelta^j_\omega, \\
B_{\omega\concat j^\infty \letterext}:=\ & B_\omega + \ZZ\cdot \diffdelta^j_\omega.
\end{align}
Here, plus denotes a Minkowski sum.
Thus for $\omega\in[n]^*$, $B_\omega$ can be regarded as a discrete zonotope, that is, an affine image of the vertices of an $|\omega|$-dimensional cube.
See Figures~\ref{figB2D} and~\ref{figB3D}.
In Section~\ref{secProofOfMainTheorem}, several equivalent descriptions of $B_\omega$ are given.
Define $Z_\omega:=-B_\omega$, and for $\filtind\in\ZZ^n$,
\begin{equation}
\label{eqSofPomega}
S(\filtind;\omega):=S^{\filtind+Z_\omega}_{\filtind+B_\omega}.
\end{equation}
Below we show that this is indeed a well-defined $S$-term (only up to natural isomorphism of course) by constructing $(b,p,q,z)\in I_4$ such that $S^{pz}_{bq}$ represents $S(\filtind;\omega)$.

Define lattices $\lattice_\omega\subseteq\ZZ^n$ for $\omega\in \letters^*_a$ inductively as follows.
Put $\lattice_\eps:=\{0\}$, $\lattice_{\omega\concat j}:=\lattice_\omega$, $\lattice_{\omega\concat j^\infty}:=\lattice_\omega$, and
\[
\lattice_{\omega\concat j^\infty \letterext}:=\lattice_\omega + \ZZ\cdot \diffdelta^j_\omega.
\]
Equivalently, $\lattice_\omega=B_\omega\cap Z_\omega = \ZZ\{\delta_\omega^j\st j\in X(\omega)\}$. 
For $\filtind,\filtind'\in\ZZ^n$ with $\filtind-\filtind'\in \lattice_\omega$, $S(\filtind;\omega)=S(\filtind';\omega)$.
Thus we may also think of $S(\filtind;\omega)$ as being parametrized over $\filtind\in\ZZ^n/\lattice_\omega$.

\begin{definition}[$\omega$-page]
Let $\omega\in \letters^*_a$.
We call the collection of all $S(\filtind;\omega)$, $\filtind\in\ZZ^n/\lattice_\omega$, the \emph{$\omega$-page}.
\end{definition}

For $\omega=\eps$ this was called the \emph{first page} in~\cite{Mat13succSpectralSequences}, for $\omega=123\ldots n$ the \emph{second page}, and for $\omega=1^{q_1}\ldots n^{q_n}$ a \emph{generalized second page}, or the $Q$-page, where $Q=(q_1,\ldots,q_n)\in\ZZ^n_{\geq 0}$.

\begin{remark}[Relation between $jj^\infty$ and $j^\infty$]
\label{remRelationBetweenJJinfAndJinf}
Suppose $w\in \letters^*_a$ contains $j^\infty$, and let $w'$ be the same word except that $j^\infty$ is replaced by $j^kj^\infty$ for some $k\geq 1$.
Then in general, $\diffdir^i_\omega\neq \diffdir^i_{\omega'}$ and $\diffdelta^i_\omega\neq \diffdelta^i_{\omega'}$, but they always agree modulo $\lattice_\omega=\lattice_{\omega'}$.
Also $B_\omega=B_{\omega'}$ and hence $S(\filtind;\omega)=S(\filtind;\omega')$.
Moreover one can check that
the differentials in the main theorem~\ref{thmDifferentials} below are the same for $\omega$ and $\omega'$.
Thus in order to speak about the $\omega$-page it is enough to know the image of $\omega$ in the quotient semigroup $L^*/(jj^\infty\sim j^\infty)$.
\end{remark}

\subsection{Main theorem}
\label{secMainTheorem}

\begin{theorem}[Main theorem]
\label{thmDifferentials}
Let $S$ denote the higher spectral sequence of a given excisive exact couple system $E$ over $D(\ZZ^n)$.
\begin{enumerate}[label={(\alph*)}]
\item\label{itMainThmFirstPage}
The first page of $S$ at $P\in\ZZ^n$ is naturally isomorphic to
\[
S(P;\eps) \iso E^p_q
\]
for any two downsets $p,q\subset\ZZ^n$ that satisfy $p = q\,\dot\cup\, \{P\}$.

\item\label{itMainThmLimit}
The limit of $S$ for any final $\omega \in L_a^*$ and $P = \ZZ^n \in \ZZ^n/\ZZ^n$ is naturally isomorphic to 
\[
S(P;\omega) \iso E^\infty_{-\infty}.
\]
\end{enumerate}
Moreover, for any $\omega\in \letters^*_a$, $P\in\ZZ^n/V_\omega$, and $j\in[n]$ such that $\omega\concat j$ is admissible, the following properties hold.
\begin{enumerate}[resume,label={(\alph*)}]
\item\label{itMainThmDifferentials}
There are natural differentials
\begin{equation}
\label{eqDifferentialsBetweenSterms}
\ldots\tof{d} S(\filtind+\diffdir_\omega^j;\omega)\tof{d} S(\filtind;\omega)\tof{d} S(\filtind-\diffdir_\omega^j;\omega)\tof{d}\ldots.
\end{equation}
Taking homology at $S(\filtind;\omega)$ yields $S(\filtind;\omega\concat j)$.

\item \label{itMainThmDifferentialsLimit}
$S(\filtind;\omega\concat j^\infty)$ is a natural subquotient of $S(\filtind;\omega\concat j^k)$ for all $k\geq 0$.

\item \label{itMainThmExtensions}
There exists a natural $\ZZ$-filtration $(F_i)_{i\in\ZZ}$ of $S(\filtind;\omega\concat j^\infty \letterext)$,
\[
0\subseteq\ldots\subseteq F_i\subseteq F_{i+1}\subseteq\ldots\subseteq S(\filtind;\omega\concat j^\infty \letterext),
\]
such that $S(\filtind+i\cdot \diffdelta^j_\omega;\omega\concat j^\infty)\iso F_i/F_{i-1}$, for all $i\in\ZZ$.
\end{enumerate}
\end{theorem}




\begin{remark}[Multiplicative structure]
\label{remMultiplicativeStructure}
As usual, under certain assumptions on $E$ there will be a multiplicative structure.
The simplest instance is when $E$ comes via~\eqref{eqExactCoupleSystemOfFilteredChainComplex} from an $I$-filtered differential algebra $C$, whose filtration $(F_p)_I$ satisfies $F_p\cdot F_q\subseteq F_{p+q}$, where $p+q$ denotes the Minkowski sum.
Then for any $\omega\in L^*_a$ there is a natural product $S(P;\omega)\otimes S(Q;\omega)\to S(P+Q;\omega)$, which satisfies a Leibniz rule with respect to the differentials~\itref{itMainThmDifferentials}.
Furthermore they are compatible with respect to~\itref{itMainThmDifferentialsLimit} and~\itref{itMainThmExtensions} in the usual way, and for final $\omega$ it coincides with the product on $H(C)$.
In other cases, such a natural product will only exist for $\omega$-pages where $\omega$ has a certain prefix; examples for this are the higher Leray--Serre spectral sequence (Example~\ref{exHigherLeraySerreSS}) and the higher Grothendieck spectral sequence below.
For details and a more general criterion see~\cite[\S{}4.4]{Mat13succSpectralSequences}.
\end{remark}

\begin{example}[Higher Leray--Serre spectral sequence]
\label{exHigherLeraySerreSS}
%
Consider a vertical tower of Serre-fibrations $F_{i-1}\incl E_{i-1} \to E_i$ ($1\leq i \leq n)$ and a generalized homology theory $h_*$, such that each $E_i$ has the homotopy type of a CW-complex.
Define $F_n:=E_n$ as if we append the trivial fibration $F_n \tof{\textrm{id}} E_n \to \pt$ to the tower.
In~\cite[\S{}5]{Mat13succSpectralSequences}, we associated to the fibration tower a spectral system over $D(\overline\ZZ^n)$.
Obviously we can restrict it to a spectral system over $D(\ZZ^n)$, and it becomes a higher spectral sequence in the sense of the present article.
Its limit is $h_*(E_0)$.
According to \cite[Thm~5.1]{Mat13succSpectralSequences}, its ``second page'' $S(P,123\ldots n)$ ($P=(p_1,\ldots,p_n)\in\ZZ^n$) can be naturally identified with
\begin{equation}
\label{eqHLSSS_2page}
S(P,123\ldots n) \iso H_{p_n}(F_n; H_{p_{n-1}}(F_{n-1};\ldots H_{p_1}(F_1; h_*(F_0))))),
\end{equation}
where the coefficients for the homology of each $F_i$ are the natural local coefficient systems over $\pi_1(F_i)$.

There is an obvious cohomological version of this higher Leray--Serre spectral sequence for the same fibration tower and for any generalized cohomology theory $h^*$.
If $h^*$ is moreover multiplicative, we obtain a product structure $S(P;\omega)\otimes S(Q;\omega)\to S(P+Q;\omega)$ for every $\omega\in L^*_a$, where it has the properties mentioned in Remark~\ref{remMultiplicativeStructure}.
This product structure agrees with the natural cup product 
of $H^*(F_n;\ldots H^*(F_1;h^*(F_0)))))$ for $\omega=123\ldots n$, and it is natural with respect to maps of fibration towers (up to homotopies) for the $\omega\in L^*_a$ with prefix $123\ldots n$.
The proof for this multiplicative structure is analogous to the arguments in~\cite[\S{}5.3]{Mat13succSpectralSequences} once we have  the description of the $\omega$-pages $S(P;\omega)$ 
as $S^{p_\omega z_\omega}_{b_\omega q_\omega}$ with $\sigma$-lexicographic downsets $z_\omega,q_\omega,p_\omega,b_\omega$ (see the proof of Theorem~\ref{thmDifferentials}).
\end{example}

\begin{example}[Higher Eilenberg--Moore spectral sequence]
\label{exHigherEilenbergMooreSS}
Consider $n$ fibrations $f_i: E_i\to B$ ($1\leq i\leq n$) over the same base space with $\pi_1(B)=0$ (as usual with Eilenberg--Moore spectral sequences, this assumption can be weakened).
Let $X$ denote the pullback of $f_1,\ldots,f_n$.
In~\cite[\S{}8]{Mat13succSpectralSequences} we constructed a spectral system over $D(\overline\ZZ^n)$ from the cube of pullbacks of $f_1,\ldots,f_n$, which becomes a higher spectral sequence in the sense of the present article.
Its limit is $H_*(X)$.
For any $\sigma\in\Sym_n$, let $\omega_\sigma = \sigma_1\concat\ldots\concat\sigma_n$.
Then the $\omega_\sigma$-page at $P=(p_1,\ldots,p_n)$ can be naturally identified with
\[
S(P;\omega_\sigma) \iso 
\Tor_{H_*(B)}^{p_{\sigma_n}}\Big(\ldots
\Tor_{H_*(B)}^{p_{\sigma_2}}\big(
\Tor_{H_*(B)}^{p_{\sigma_1}}\big(
H_*(B),H_*(E_{\sigma_1})\big)
,H_*(E_{\sigma_2})\big)
\ldots,H_*(E_{\sigma_n})\Big).
\]
Here, $\Tor_{H_*(B)}^{p_{\sigma_1}}(
H_*(B),H_*(E_{\sigma_1}))$ simplifies to $H_*(E_{\sigma_1})$ if $p_{\sigma_1}=0$ and to $0$ otherwise.
Note that in this example, we have $n!$ meaningful ``second pages'', in comparison to just one in Example~\ref{exHigherLeraySerreSS}.
\end{example}

\begin{example}[Higher Grothendieck spectral sequence]
Let 
$\cat A_0 \tof{F_1} \cat A_1 \tof{F_2} \ldots \tof{F_n} \cat A_n$
be a sequence of $n$ right exact functors between Abelian categories with enough projectives.
Write $F_{ij} = F_j \circ \ldots \circ F_{i+1}: \cat A_i \to \cat A_j$ and assume that for any $i<j<k$, $F_{ij}$ sends projective objects to $F_{jk}$-acyclic objects.
In~\cite[\S{}6]{Mat13succSpectralSequences} we constructed for each object $X\in\cat A_0$ a spectral system over $D(\overline\ZZ^n)$, which becomes a higher spectral sequence in the sense of the present article. Its limit is $L_*F_{0n}(X)$.
Its ``second page'' $S(P;123\ldots n)$ ($P=(p_1,\ldots,p_n)\in\ZZ^n)$) can be naturally identified with
\[
S(P;123\ldots n) \iso
(L_{p_n}F_n) \circ \ldots \circ (L_{p_1}F_1)(X).
\]
Furthermore, in the setting of~\cite[\S{}6.3]{Mat13succSpectralSequences} this higher Grothendieck spectral sequence admits a natural product structure $S(P;\omega)\otimes S(Q;\omega)\to S(P+Q;\omega)$ for every $\omega\in L^*_a$ with prefix $123\ldots n$; and for such $\omega$ it has the properties mentioned in Remark~\ref{remMultiplicativeStructure}.
\end{example}

Similarly, the spectral system obtained from 
the chromatic spectral sequence followed by the Adams--Novikov spectral sequence~\cite[\S{}7]{Mat13succSpectralSequences} becomes a higher spectral sequence in the sense of the present article.

\begin{example}[Universal coefficients for spectral sequences]
Let $C$ be a free graded chain complex over a ring~$R$ and let $M$ be an $R$-module.
There is a universal coefficient spectral sequence that relates $H_*(C)$ with $H_*(C\otimes_R M)$, which simplifies to the usual universal coefficient theorem for $\blank\otimes_R M$ in case $R$ is a principal ideal domain.
This motivates the question, whether there exists a universal coefficient theorem that relates the spectral sequence of a $\ZZ$-filtration $(F_n)_n$ of $C$ with the spectral sequence of the induced filtration $(F_n\otimes_R M)_n$ of $C\otimes_R M$.
Let us assume that the subchain complexes $F_n$ are free.
Then a reasonable positive answer to this is the following ``universal coefficient higher spectral sequence for $\blank\otimes_R M$''.
Let $P_*\to M$ be a projective resolution of $M$.
Define the (total) chain complex $D:=C\otimes_R P_*$, which is filtered in two ways, $(F_n\otimes_R P_*)_n$ and $(C\otimes_R P_{\leq m})_m$.
They refine to an $I$-filtration ($I=D(\ZZ^2)$), which is distributive and hence the associated exact couple system is excisive.
The induced higher spectral sequence has limit naturally isomorphic to $H_*(C\otimes_R M)$.
Various of its $\omega$-pages can be naturally identified, such as $S((n,m);1^{r-1}2) = \Tor_m(E^r_n(C),M)$ ($1\leq r\leq \infty$) and $S([(0,m)];1^\infty\letterext2) = \Tor_m(H_*(C),M)$, where $E^r_n(C)$ denotes the $n$'th column of the $r$'th page of the spectral sequence of $C$ with filtration $(F_n)_n$.
Analogously there is a ``universal coefficient higher spectral sequence for $\Hom_R(\blank,M)$''.

More generally, every higher spectral sequence of a graded chain complex $C$ with distributive $D(\ZZ^n)$-filtration $(F_p)_p$ yields analogously two universal coefficient higher spectral sequences over $D(\ZZ^{n+1})$ that relate to the higher spectral sequences over $D(\ZZ^n)$ of $(F_p\otimes_R M)_p$ and $(\Hom_R(F_p,M))_p$, respectively.
Another or additional way of generalizing is to replace $M$ in the above by another free graded chain complex, possibly with a $\ZZ$- or $I$-filtration of its own, and we arrive at various higher spectral sequences of K{\"u}nneth type.
\end{example}

\subsection{Proof of the main theorem}
\label{secProofOfMainTheorem}

Before proving Theorem~\ref{thmDifferentials} we need some preparation.
In particular it will be convenient to move to another basis spanned by the $\diffdelta^i_\omega$, which depends on $\omega$.

For $\omega\in \letters^*_a$, define
\[
\transf_\omega:=\begin{pmatrix}\diffdelta^1_\omega & \cdots & \diffdelta^n_\omega\end{pmatrix}^{-1}\in \QQ^{n\times n}.
\]
For $k\in[n]$ and $\omega\in\letters^*_a$, let $\transfdelta^k_\omega$ denote the matrix
\[
\transfdelta^k_\omega:=\id_{\ZZ^n}+e_k(\one_{Y(\omega)\wo\{k\}})^\tr
= (\delta_{i=j\textnormal{ or }(i=k\textnormal{ and }j\in Y(\omega))})_{ij}
\in \SL(n,\ZZ).
\]
A quick calculation shows that $\transf_\eps = \id_n$, $\transf_{\omega\concat j^\infty}=\transf_\omega$, $\transf_{\omega\concat \letterext}=\transf_\omega$, and
\[
\transf_{\omega\concat j} = \transfdelta^j_\omega\transf_\omega.
\]
Therefore, $\transf_\omega\in\SL(n,\ZZ)$, and all its entries are non-negative.
It will be very convenient to transform $B_\omega$ and $\lattice_\omega$ by $\transf_\omega$, so we define
\[
B^\transf_\omega:=\transf_\omega\cdot B_\omega,\ \  \lattice^\transf_\omega:=\transf_\omega\cdot \lattice_\omega\ \ \subseteq \ZZ^n.
\]
In particular,
\begin{align}
\label{eqBTinduction}
B^\transf_{\omega\concat j} =\ & \transfdelta^j_\omega(B^\transf_\omega + \{0,e_j\}), \\
B^\transf_{\omega\concat j^\infty} =\ & \transfdelta^j_\omega(B^\transf_\omega + \ZZ_{\geq 0}\cdot e_j), \\
B^\transf_{\omega\concat j^\infty \letterext} =\ & \transfdelta^j_\omega(B^\transf_\omega + \ZZ\cdot e_j),
\end{align}
and $\lattice^\transf_\omega = \ZZ\cdot\{e_i\st i\in X(\omega)\}$.
In particular, $\ZZ^n/\lattice^\transf_\omega$ can be naturally identified with~$\ZZ^{Y(\omega)}$.



Define $u_\omega\in\ZZ^n_{\geq 0}$ inductively via $u_\eps:=0$, $u_{\omega\concat j^\infty}:=u_\omega$, $u_{\omega\concat\letterext}:=u_\omega$, and
\[
u_{\omega\concat j}:=e_j+\transfdelta^j_\omega u_\omega.
\]
For $n\geq 2$ and $\omega\in[n]^*$, a more explicit formula is $\upperbound_\omega=\tfrac{1}{n-1}(\transf_\omega\one-\one)$.

The following lemmas about $B_\omega$ are stated and proved only in the special case when $\omega\in[n]^*$, since this simplifies the notation considerably.
When working modulo $\lattice_\omega$ and $\lattice^\transf_\omega$, respectively, analogous statements still hold for arbitrary $\omega\in \letters^*_a$, as long as $\omega$ does not end on some $j^\infty$. 
It should be rather clear how to state and prove them.

\begin{lemma}
For any $\omega\in [n]^*$, $B^\transf_\omega$ is a path in the unit-distance graph of $\ZZ^n$ from $0$ to $\upperbound_\omega$, which is monotone with respect to all coordinates.
\end{lemma}

We remark that this path is in general not simply a discretized line segment (i.e., the set of all lattice points with $\ell^\infty$-distance at most $1/2$ from some line segment in~$\RR^n$).

\begin{proof}
We proceed by induction. 
Assume that the lemma holds for $B^\transf_\omega$, and we want to prove it for~$B^\transf_{\omega\concat j}$ using~\eqref{eqBTinduction}.
Suppose $x,x+e_k\in B^\transf_\omega$ are the vertices of an edge in $B^\transf_\omega$.
If $k\neq j$ then this edge gives rise to a path of length 2 in $B^\transf_{\omega\concat j}$ along the vertices $\transfdelta^j_\omega x$, $\transfdelta^j_\omega(x+e_j)=\transfdelta^j_\omega x+e_j$, and $\transfdelta^j_\omega(x+e_k)=\transfdelta^j_\omega x+e_j+e_k$.
If $k=j$ then this edge gives rise to an edge in $B^\transf_{\omega\concat j}$ whose vertices are $\transfdelta^j_\omega x$ and $\transfdelta^j_\omega(x+e_k)=\transfdelta^j_\omega x+e_j$.
Moreover $\upperbound_{\omega\concat j} = \transfdelta^j_\omega(\upperbound_\omega+e_j)$.
\end{proof}

\begin{corollary}
\label{corSymmetryOfB}
For any $\omega\in [n]^*$, $B^\transf_\omega\subseteq \{x\in\ZZ^n\st 0\leq x\leq \upperbound_\omega\}$.
Moreover, $B^\transf_\omega = \upperbound_\omega-B^\transf_\omega$.
\end{corollary}

\begin{lemma}
\label{lemNeighborhoodOfB}
For any $\omega\in [n]^*$, $i,j\in[n]$, $x\in B^\transf_\omega$, the following holds:
\begin{enumerate}
\item Either $x-\transf_\omega e_i\in B^\transf_\omega$, or $x-\transf_\omega e_i\leq 0$, or both.
\item Either $x+\transf_\omega e_i\in B^\transf_\omega$, or $x+\transf_\omega e_i\geq \upperbound_\omega$, or both.
\end{enumerate}
\end{lemma}

\begin{proof}
Suppose the lemma holds for $B^\transf_\omega$, and we want to prove it for $B^\transf_{\omega\concat j}$.
By Corollary~\ref{corSymmetryOfB} it is enough to prove the first statement.
By~\eqref{eqBTinduction}, any element of $B^\transf_{\omega\concat j}$ is of the form $\transfdelta^j_\omega x$ or $\transfdelta^j_\omega(x+e_j)=\transfdelta^j_\omega x+e_j$ for some $x\in B^\transf_\omega$.
If now $z:=x-\transf_\omega e_i\in B^\transf_\omega$, then also $z':=\transfdelta^j_\omega x-\transf_{\omega\concat j}e_i=\transfdelta^j_\omega z$ and $z'':=\transfdelta^j_\omega(x+e_j)-\transf_{\omega\concat j}e_i=\transfdelta^j_\omega(z+e_j)$ lie in $B^\transf_{\omega\concat j}$.

Thus it remains to check the case when $z\leq 0$.
Then clearly $z'=\transfdelta^j_\omega z\leq 0$ since all entries of $\transfdelta^j_\omega$ are non-negative.
Similarly $z''=\transfdelta^j_\omega z+e_j\leq e_j$.
If $z''\leq 0$ does not hold, then $1=z''_j=1+\sum_kz_k$, hence $z=0$, thus $z''=e_j$, which lies in $B^\transf_{\omega\concat j}$ since $0\in B^\transf_\omega$.
\end{proof}

Let us regard any subset of $\ZZ^n$ as a graph by connecting any two elements with distance $1$ by an edge.
In particular we can then talk about connected components of such subsets.

\begin{lemma}
\label{lemBisConnected}
For any $\omega\in [n]^*$, $B_\omega$ is connected.
\end{lemma}

\begin{proof}
By induction one immediately sees that
\begin{equation}
\label{eqAjMinusEjIsInB}
\diffdir^j_\omega-e_j\in B_\omega.
\end{equation}
Also by induction,
\begin{equation}
\label{eqTAisUplusEj}
\transf_\omega \diffdir_\omega^j = \upperbound_\omega+e_j\in e_j+B^\transf_\omega,
\end{equation}
from which we get
\begin{equation}
\label{eqAjIsInRjPlusB}
\diffdir_\omega^j\in \diffdelta^j_\omega+B_\omega.
\end{equation}
Thus if $B_\omega$ is connected then so is $B_{\omega\concat j}$ by~\eqref{eqBasZonotope}, \eqref{eqAjMinusEjIsInB}, and~\eqref{eqAjIsInRjPlusB}.
\end{proof}

For any subset $X\subseteq \ZZ^n$, let $\Comp_0(X)$ denote the connected component that contains~$0$.

\begin{corollary}
For any $\omega\in[n]^*$, $B_\omega$ is $\Comp_0(X)$ where $X$ is the intersection of the $n$ ``discrete hyperplanes''
\[
\{x\in\ZZ^n\st 0\leq e_i^\tr \transf_\omega x \leq e_i^\tr \upperbound_\omega\},\quad 1\leq i\leq n.
\]
\end{corollary}

\begin{proof}
This follows from Lemmas~\ref{lemNeighborhoodOfB} and~\ref{lemBisConnected}.
\end{proof}

For our purposes the following similar description of $B_\omega$ is more useful.
The symmetric group $\Sym_n$ acts on $\ZZ^n$ by permutation of the coordinates, $\sigma\cdot x:=(x_{\sigma^{-1}(i)})_{i\in[n]}$.
Let $\leqlex$ denote the lexicographic relation on $\ZZ^n$.
For $\sigma\in\Sym_n$, define a new relation $\leqlexsigma{\sigma}$ by setting $x\leqlexsigma{\sigma} y$ if and only if $\sigma x\leqlex \sigma y$.

For any $\filtind\in (\ZZ\cup\{\infty\})^n$, $\transf\in\ZZ^{n\times n}_{\geq 0}$ with $\det \transf=1$, and $\sigma\in\Sym_n$, define
\[
\downset(\filtind;\transf,\sigma):=\{x\in\ZZ^n\st \transf x\leqlexsigma{\sigma} \transf \filtind\}
\]
and
\[
\downset^\circ(\filtind;\transf,\sigma):=\{x\in\ZZ^n\st \transf x\lesslexsigma{\sigma} \transf \filtind\}.
\]

\begin{lemma}
\label{lemBviaLexicographicDownsets}
For any $\omega\in[n]^*$, $j\in[n]$, and $\sigma\in\Sym_n$ with $\sigma(j)=n$, the following two equations hold.
\begin{enumerate}
\item 
$
B_\omega = \Comp_0\big(\downset^\circ(\diffdir^j_\omega; \transf_\omega,\sigma)\wo \downset^\circ(0; \transf_\omega,\sigma) \big).
$
\item 
$
B_{\omega\concat j} = \Comp_0\big(\downset(\diffdir^j_\omega; \transf_\omega,\sigma)\wo \downset^\circ(0; \transf_\omega,\sigma) \big).
$
\end{enumerate}
\end{lemma}

Of course analogous formulas also hold for $Z_\omega$ and $Z_{\omega\concat j}$.

\begin{proof}
Since $\sigma(j)=n$, the only $z\in\ZZ^n$ with $\upperbound_\omega\lesslexsigma{\sigma} \transf_\omega z \leqlexsigma{\sigma} \transf_\omega \diffdir^j_\omega$ is $\diffdir^j_\omega$.
Now the first equation follows readily from Lemmas~\ref{lemNeighborhoodOfB} and~\ref{lemBisConnected}.

As for the second equation, ``$\subseteq$'' follows from $\transf_\omega B_{\omega\concat j} = \transf_\omega B_\omega + \{0,e_j\}$ and the connectivity of $B_{\omega\concat j}$.
It remains to check ``$\supseteq$''.
By Lemma~\ref{lemNeighborhoodOfB} and~\eqref{eqTAisUplusEj} the only neighbor~$y$ of $x\in B_\omega$ such that $y\not\in B_\omega$ and $y\in \downset(\diffdir^j_\omega; \transf_\omega,\sigma)\wo \downset^\circ(0; \transf_\omega,\sigma)$ is $y=\diffdir^j_\omega$.
Similarly (or by symmetry), the only neighbor~$y$ of $x\in B_\omega+\diffdelta^j_\omega$ such that $y\not\in B_\omega+\diffdelta^j_\omega$ and $y\in \downset(\diffdir^j_\omega; \transf_\omega,\sigma)\wo \downset^\circ(0; \transf_\omega,\sigma)$ is $y=0$.
Both, $0$ and $\diffdir^j_\omega$, lie in $B_{\omega\concat j}$, which proves the claimed equality.
\end{proof}

\begin{proof}[Proof of Theorem~\ref{thmDifferentials}]
We may assume $\filtind=0$, otherwise translate everything.
\itref{itMainThmFirstPage} and \itref{itMainThmLimit} hold by definition and the excision property.

\itref{itMainThmDifferentials}
We first consider the case $\omega\in[n]^*$.
Let $\sigma\in\Sym_n$ be any permutation with $\sigma(j)=n$.
Define 
\begin{align*}
p_\omega :=\ & \downset(0; \transf_\omega,\sigma), \\
q_\omega :=\ & \downset^\circ(0; \transf_\omega,\sigma) = p_\omega\wo 0, \\
b_\omega :=\ & \downset^\circ(\diffdir^j_\omega;\transf_\omega,\sigma), \\
z_\omega :=\ & \downset(-\diffdir^j_\omega;\transf_\omega,\sigma).
\end{align*}
Then Lemma~\ref{lemBviaLexicographicDownsets} implies that $S^{p_\omega z_\omega}_{b_\omega q_\omega}$ represents $S^{Z_\omega}_{B_\omega}$, which by definition is $S(0;\omega)$.
One can describe $S(\diffdir^j_\omega;\omega)$ and $S(-\diffdir^j_\omega;\omega)$ similarly by translating all downsets by $\pm \diffdir^j_\omega$.
Now Lemma~\ref{lemDifferentials} implies the claimed differentials in~\eqref{eqDifferentialsBetweenSterms}.
Lemma~\ref{lemKernelsAndCokernels} shows that taking homology in~\eqref{eqDifferentialsBetweenSterms} at $S(0;\omega)$ yields $S^{p_\omega z_\omega^*}_{b_\omega^* q_\omega}$ where 
\begin{align*}
b_\omega^* :=\ & \downset(\diffdir^j_\omega;\transf_\omega,\sigma) = b_\omega\cup\{\diffdir^j_\omega\}, \\
z_\omega^* :=\ & \downset^\circ(-\diffdir^j_\omega;\transf_\omega,\sigma) = z_\omega\wo\{-\diffdir^j_\omega\}.
\end{align*}
By Lemma~\ref{lemBviaLexicographicDownsets}, $S^{p_\omega z_\omega^*}_{b_\omega^* q_\omega}$ represents $S^{Z_{\omega\concat j}}_{B_{\omega\concat j}} = S(0;\omega\concat j)$.


The general case $\omega\in \letters^*_a$ only needs minor modifications:
In particular, the downsets $b_\omega$, $p_\omega$, $q_\omega$, $z_\omega$ need to be replaced by the sum of $\lattice_\omega$ with the analogous downsets in $\ZZ^n/\lattice_\omega$.
Explicitly, let
$k:=|X(\omega)|$, and choose a $\sigma\in\Sym_n$ that satisfies $\sigma(j)=n-k$ and $\sigma(i)>n-k$ for all $i\in X(\omega)$.
Then $p_\omega$ can be defined as 
\[
p_\omega := \{x\in\ZZ^n\st \transf_\omega x\leqlexsigma{\sigma} (0^{n-k},\infty^{k})\},
\]
and put $q_\omega:=p_\omega\wo \lattice_\omega$, $b_\omega:=q_\omega+\diffdir^j_\omega$, $z_\omega:=p_\omega-\diffdir^j_\omega$, $b_\omega^*:=p_\omega+\diffdir^j_\omega=b_\omega\cup(\diffdir^j_\omega+\lattice_\omega)$ and $z_\omega^*:=q_\omega-\diffdir^j_\omega=z_\omega\wo(-\diffdir^j_\omega+\lattice_\omega)$.
Note that $\lattice^\transf_\omega = \ZZ^{X(\omega)}$.
Now one can repeat the previous argument in the quotient space $\ZZ^n/\lattice_\omega$.

\smallskip

\itref{itMainThmDifferentialsLimit}
In $\transf_\omega$-coordinates,
\[
\transf_\omega B_{\omega\concat j^i} = B^\transf_\omega + \{0,\ldots,i\}.
\]
Hence as in the proof of Lemma~\ref{lemBviaLexicographicDownsets} one can show that 
\[
B_{\omega\concat j^i} = \Comp_0\big((b_\omega + i\cdot \diffdelta_\omega^i)\wo q_\omega\big),
\]
where $b_\omega$ and $q_\omega$ are as above.
Therefore $S(0;\omega\concat j^i)$ can be represented by $S^{p_\omega z_\omega^i}_{b_\omega^i q_\omega}$, where
\begin{align*}
b_\omega^i :=\ & b_\omega+i\cdot \diffdelta^j_\omega 
     = \{x\in\ZZ^n\st \transf_\omega x\leqlexsigma{\sigma} \upperbound_\omega+(0^{n-k-1},i,\infty^{k})^\tr\}, \\
z_\omega^i :=\ & z_\omega-i\cdot \diffdelta^j_\omega 
     = \{x\in\ZZ^n\st \transf_\omega x\leqlexsigma{\sigma} -\upperbound_\omega+(0^{n-k-1},-i-1,\infty^{k})^\tr\}. 
\end{align*}
Let $b_\omega^\infty:=\bigcup_i b_\omega^i$ and $z_\omega^\infty:=\bigcap_i z_\omega^i$.
Then $S(0;\omega\concat j^\infty) = S^{p_\omega z_\omega^\infty}_{b_\omega^\infty,q_\omega}$.

\smallskip

\itref{itMainThmExtensions}
Let $p_\omega,b_\omega^\infty,z_\omega^\infty$ as above.
Put $p_\omega^i:=i\cdot \diffdelta^j_\omega + p_\omega$, $p_\omega^{-\infty}:=\bigcap_{i\in\ZZ} p_\omega^i$, and $p_\omega^{\infty}:=\bigcup_{i\in\ZZ} p_\omega^i$.
Define
\[
F_i:=\im\Big(\ell: S^{p_\omega^i z_\omega^\infty}_{b_\omega^\infty p_\omega^{-\infty}}\to S^{p_\omega^\infty z_\omega^\infty}_{b_\omega^\infty p_\omega^{-\infty}}\Big).
\]
Then the assertion follows from Lemma~\ref{lemExtensionProperty} (or from the proof of Lemma~\ref{lemExactCouplesInfinityPageFiltrationQuotients}).
\end{proof}

\subsection{The $2$-dimensional case}
\label{secCaseN2}


The probably most frequent case (apart from the classical one, $n=1$) is $n=2$.
A few more things can be said about this case:


Every final $\omega\in \letters^*_a$ is of the form 
\begin{equation}
\label{eqFormOfOmegaForNequalTo2}
\omega=\tau\concat j_1^\infty\letterext j_2^kj_2^\infty\letterext,
\end{equation}
for some $\tau\in[2]^*$, $\{j_1,j_2\}=[2]$, and $k\geq 0$.
Any such $\omega$ gives a recipe to connect the first page to the limit of the spectral system.
This recipe is therefore already determined by $\tau$ and~$j_1$.

Note that for all prefixes $\tau\leq\omega'\leq\omega$, $M_\tau = M_{\omega'}=M_\omega$.
Let's define
\[
N_\omega := (e_{j_2})^\tr\cdot M_\omega,
\]
which is the ``normal vector'' along which the downsets $b,p,q,z$ grow respectively shrink during $j_2^kj_2^\infty\letterext$.
Clearly $N_\omega\geq 0$ and it is primitive (i.e.\ its entries are coprime), and $N=(e_i)^\tr$ can happen only if $i=j_2$.
Also, $N_\omega$ is invariant under the relation $jj^\infty\sim j^\infty$, compare with Remark~\ref{remRelationBetweenJJinfAndJinf}.

\begin{observation}
Modulo $jj^\infty\sim j^\infty$, $\omega$ is uniquely determined by $N_\omega$ and $j_1$.
Conversely, for any primitive $N^\tr\in\ZZ^2_{\geq 0}$ and $j_1\in[2]$ with $N^\tr\neq e_{j_1}$ there is a final $\omega\in L^*_a$ of the form~\eqref{eqFormOfOmegaForNequalTo2} such that $N=N_\omega$.
\end{observation}

Thus the connection determined by $\omega$ can be equivalently described by the pair $(N_\omega,j_1)$.

\begin{proof}
In fact there is a simple algorithm that determines all possible $\tau$ from $N$ (respectively $N_\omega$) and $j_2=3-j_1$.
If $j_2=1$, choose $(N')^\tr\in\ZZ^2_{\geq 0}$ such that
$M:=\bigl(\begin{smallmatrix}
N\\ N'
\end{smallmatrix} \bigr) \in \SL(2,\ZZ)\cap \ZZ_{\geq 0}^{2\times 2}$.
If $j_2=2$, choose $(N')^\tr\in\ZZ^2_{\geq 0}$ such that
$M:=\bigl(\begin{smallmatrix}
N'\\ N
\end{smallmatrix} \bigr) \in \SL(2,\ZZ)\cap \ZZ_{\geq 0}^{2\times 2}$.
In any case, $N'$ is well-defined up to adding an integral multiple of~$N$.
If $N'_0$ is the smallest choice, then all others are of the form $N'_k:=N'_0+kN$, $k\in\ZZ_{\geq 0}$.
Now one can repetitively take one of the two rows of $M$ and subtract it from the other one such that all entries stay non-negative until one arrives at $\id_{\ZZ^2}$, and there is a unique way to do that.
Let $q_i\in[2]$ denote the index of the column from which the other column was subtracted during the $i$'th round. And say there were $\ell$ rounds.
Then $M = M_\omega$ for $\omega=\tau j_1^\infty\letterext j_2^\infty\letterext$ and $\tau:=q_\ell\concat\ldots\concat q_1$.
The choice of $N'$ correspond to how often $j_1$ appears at the end of~$\tau$, namely $k$ times if $N'=N'_k$.
\end{proof}

The algorithm has similarities to the extended Euclidean algorithm applied to the first column of~$M$.
The only difference is that in the extended Euclidean algorithm, one subtracts \emph{multiples} of one number from the other.

\begin{example}
Consider an excisive exact couple system $E$ over $I(\ZZ^2)$.
Suppose we want to determine both $\omega$ for which $N_\omega=(3,5)$.
For $j_1=1$, the algorithm runs as follows:
\[
\bigl(\begin{smallmatrix}
N' \\ N
\end{smallmatrix} \bigr)
=
\bigl(\begin{smallmatrix}
2 & 3 \\ 3 & 5
\end{smallmatrix} \bigr)
\Tof{q_1=2}
\bigl(\begin{smallmatrix}
2 & 3 \\ 1 & 2
\end{smallmatrix} \bigr)
\Tof{q_2=1}
\bigl(\begin{smallmatrix}
1 & 1 \\ 1 & 2
\end{smallmatrix} \bigr)
\Tof{q_3=2}
\bigl(\begin{smallmatrix}
1 & 1 \\ 0 & 1
\end{smallmatrix} \bigr)
\Tof{q_4=1}
\bigl(\begin{smallmatrix}
1 & 0 \\ 0 & 1
\end{smallmatrix} \bigr),
\]
Thus $\omega=12121^\infty\letterext 2^\infty\letterext$ does it.
Similarly, for $j_1=2$,
\[
\bigl(\begin{smallmatrix}
N \\ N'
\end{smallmatrix} \bigr)
=
\bigl(\begin{smallmatrix}
3 & 5 \\ 1 & 2
\end{smallmatrix} \bigr)
\Tof{q_1=1}
\bigl(\begin{smallmatrix}
2 & 3 \\ 1 & 2
\end{smallmatrix} \bigr)
\Tof{q_2=1}
\bigl(\begin{smallmatrix}
1 & 1 \\ 1 & 2
\end{smallmatrix} \bigr)
\Tof{q_3=2}
\bigl(\begin{smallmatrix}
1 & 1 \\ 0 & 1
\end{smallmatrix} \bigr)
\Tof{q_4=1}
\bigl(\begin{smallmatrix}
1 & 0 \\ 0 & 1
\end{smallmatrix} \bigr),
\]
and thus $\omega=12112^\infty\letterext 1^\infty\letterext$ does it too. See Figure~\ref{figN35}.
\begin{figure}[htb]
\centerline{
\input{exampleNis35v4.pspdftex}
} 
\caption{
$B_{12121^\infty}$, $B_{12112^\infty}$, and $B_{12121^\infty\letterext 2^k}=B_{12112^\infty\letterext 1^k}$ for $0\leq k\leq 5$ and for $k=\infty$.
In the third figure, the squares with number $i$ belong to $B_{12121^\infty\letterext 2^k}$ if and only if $i\leq k$.
The solid squares depict~$V_\sigma$.}
\label{figN35}
\end{figure}
\end{example}


\paragraph{Continued fractions.}
This algorithm is essentially the same as the one behind continued fractions:
Write
\[
(a_0,a_1,\ldots,a_\ell) := a_0 +\frac{1}{a_1+\frac{1}{\ldots+\frac{1}{a_\ell}}}.
\]
Let $N^\tr=(x,y)\in \ZZ^2_{\geq 0}$ be primitive.
Write the slope of $N$ as a continued fraction, $\frac{y}{x}=(a_0,\ldots,a_\ell,\infty)$, with $a_i\in\ZZ$, all positive except for possibly $a_0=0$.
For each $N\neq e_i^\tr$, there are two such representations, namely one with $a_\ell\geq 2$, and one with $a_\ell=1$: To get from the former representation to the latter, write $a_\ell = (a_\ell-1)+\frac{1}{1}$.
Comparing the above algorithm with the recursion for the successive convergents of this continued fraction, one sees immediately that $N=N_\omega$ for $\tau = 1^{a_0}2^{a_1}1^{a_2}2^{a_3}\cdots$.

\begin{funfact}
The golden ratio can be arbitrarily well approximated by the slope of $N_\omega$ using $\tau=(12)^k$, since for this $\tau$,
$M_\tau=
\bigl(\begin{smallmatrix}
1 & 1 \\ 1 & 2
\end{smallmatrix} \bigr)^k =
\bigl(\begin{smallmatrix}
f_{2k-1} & f_{2k} \\ f_{2k} & f_{2k+1}
\end{smallmatrix} \bigr)$,
where $f_k$ are the Fibonacci numbers.
In terms of continued fractions, this is because the golden ratio satisfies $(1+\sqrt{5})/2 = (1,1,1,\ldots)$.
However irrational slopes are not particularly useful, since one cannot connect the obtained page naturally to the limit (at least without further assumptions on $E$ and without going backwards).
\end{funfact}


%
%
%


%
%

\small

\bibliographystyle{plain}


\bibliography{../mybib07}

\end{document}

%% file: differentialTree2D3.pspdftex
\begin{picture}(0,0)%
\includegraphics{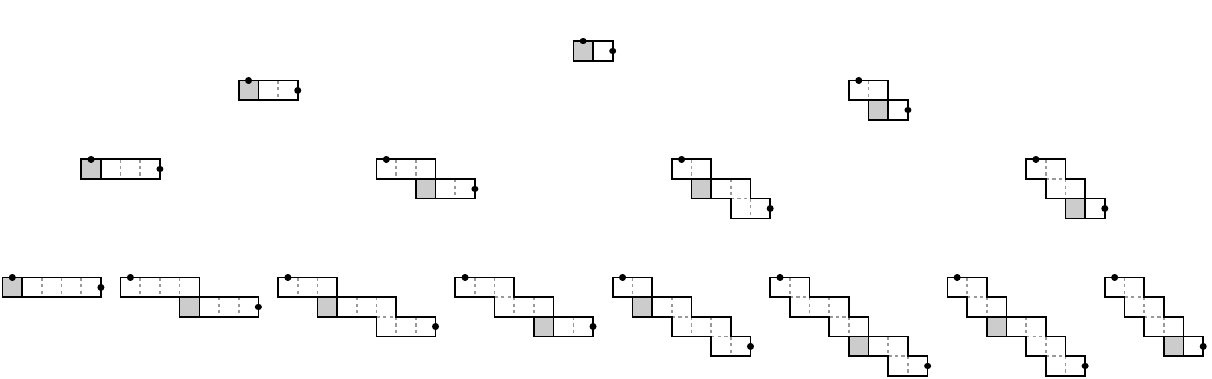}%
\end{picture}%
\setlength{\unitlength}{4144sp}%
\begingroup\makeatletter\ifx\SetFigFont\undefined%
\gdef\SetFigFont#1#2#3#4#5{%
  \reset@font\fontsize{#1}{#2pt}%
  \fontfamily{#3}\fontseries{#4}\fontshape{#5}%
  \selectfont}%
\fi\endgroup%
\begin{picture}(5535,1719)(15379,-13933)
\put(20521,-13381){\makebox(0,0)[lb]{\smash{{\SetFigFont{6}{7.2}{\rmdefault}{\mddefault}{\updefault}{\color[rgb]{0,0,0}$1222$}%
}}}}
\put(19801,-13381){\makebox(0,0)[lb]{\smash{{\SetFigFont{6}{7.2}{\rmdefault}{\mddefault}{\updefault}{\color[rgb]{0,0,0}$1221$}%
}}}}
\put(18991,-13381){\makebox(0,0)[lb]{\smash{{\SetFigFont{6}{7.2}{\rmdefault}{\mddefault}{\updefault}{\color[rgb]{0,0,0}$1212$}%
}}}}
\put(18271,-13381){\makebox(0,0)[lb]{\smash{{\SetFigFont{6}{7.2}{\rmdefault}{\mddefault}{\updefault}{\color[rgb]{0,0,0}$1211$}%
}}}}
\put(17551,-13381){\makebox(0,0)[lb]{\smash{{\SetFigFont{6}{7.2}{\rmdefault}{\mddefault}{\updefault}{\color[rgb]{0,0,0}$1122$}%
}}}}
\put(16741,-13381){\makebox(0,0)[lb]{\smash{{\SetFigFont{6}{7.2}{\rmdefault}{\mddefault}{\updefault}{\color[rgb]{0,0,0}$1121$}%
}}}}
\put(16021,-13381){\makebox(0,0)[lb]{\smash{{\SetFigFont{6}{7.2}{\rmdefault}{\mddefault}{\updefault}{\color[rgb]{0,0,0}$1112$}%
}}}}
\put(15481,-13381){\makebox(0,0)[lb]{\smash{{\SetFigFont{6}{7.2}{\rmdefault}{\mddefault}{\updefault}{\color[rgb]{0,0,0}$1111$}%
}}}}
\put(15841,-12841){\makebox(0,0)[lb]{\smash{{\SetFigFont{6}{7.2}{\rmdefault}{\mddefault}{\updefault}{\color[rgb]{0,0,0}$111$}%
}}}}
\put(17191,-12841){\makebox(0,0)[lb]{\smash{{\SetFigFont{6}{7.2}{\rmdefault}{\mddefault}{\updefault}{\color[rgb]{0,0,0}$112$}%
}}}}
\put(18541,-12841){\makebox(0,0)[lb]{\smash{{\SetFigFont{6}{7.2}{\rmdefault}{\mddefault}{\updefault}{\color[rgb]{0,0,0}$121$}%
}}}}
\put(20161,-12841){\makebox(0,0)[lb]{\smash{{\SetFigFont{6}{7.2}{\rmdefault}{\mddefault}{\updefault}{\color[rgb]{0,0,0}$122$}%
}}}}
\put(16516,-12481){\makebox(0,0)[lb]{\smash{{\SetFigFont{6}{7.2}{\rmdefault}{\mddefault}{\updefault}{\color[rgb]{0,0,0}$11$}%
}}}}
\put(19306,-12481){\makebox(0,0)[lb]{\smash{{\SetFigFont{6}{7.2}{\rmdefault}{\mddefault}{\updefault}{\color[rgb]{0,0,0}$12$}%
}}}}
\put(18046,-12301){\makebox(0,0)[lb]{\smash{{\SetFigFont{6}{7.2}{\rmdefault}{\mddefault}{\updefault}{\color[rgb]{0,0,0}$1$}%
}}}}
\end{picture}%

%% file: exampleNis35v4.pspdftex
\begin{picture}(0,0)%
\includegraphics{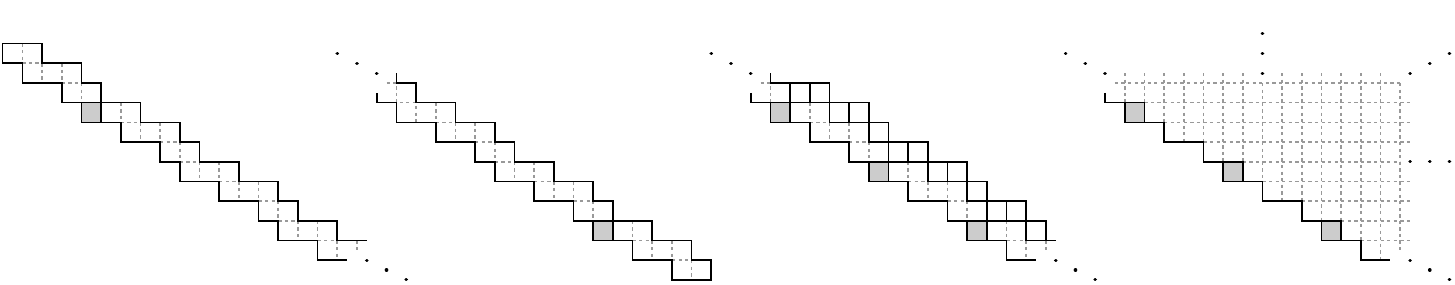}%
\end{picture}%
\setlength{\unitlength}{4144sp}%
\begingroup\makeatletter\ifx\SetFigFont\undefined%
\gdef\SetFigFont#1#2#3#4#5{%
  \reset@font\fontsize{#1}{#2pt}%
  \fontfamily{#3}\fontseries{#4}\fontshape{#5}%
  \selectfont}%
\fi\endgroup%
\begin{picture}(6649,1291)(14929,-19883)
\put(15301,-18691){\makebox(0,0)[lb]{\smash{{\SetFigFont{6}{7.2}{\rmdefault}{\mddefault}{\updefault}{\color[rgb]{0,0,0}$12121^\infty$}%
}}}}
\put(17011,-18691){\makebox(0,0)[lb]{\smash{{\SetFigFont{6}{7.2}{\rmdefault}{\mddefault}{\updefault}{\color[rgb]{0,0,0}$12112^\infty$}%
}}}}
\put(18271,-18691){\makebox(0,0)[lb]{\smash{{\SetFigFont{6}{7.2}{\rmdefault}{\mddefault}{\updefault}{\color[rgb]{0,0,0}$12121^\infty\letterext 2^k$ and $12112^\infty\letterext 1^k$}%
}}}}
\put(19891,-18691){\makebox(0,0)[lb]{\smash{{\SetFigFont{6}{7.2}{\rmdefault}{\mddefault}{\updefault}{\color[rgb]{0,0,0}$12121^\infty\letterext 2^\infty$ and $12112^\infty\letterext 1^\infty$}%
}}}}
\put(18587,-19040){\makebox(0,0)[b]{\smash{{\SetFigFont{6}{7.2}{\rmdefault}{\mddefault}{\updefault}{\color[rgb]{0,0,0}$1$}%
}}}}
\put(19037,-19310){\makebox(0,0)[b]{\smash{{\SetFigFont{6}{7.2}{\rmdefault}{\mddefault}{\updefault}{\color[rgb]{0,0,0}$1$}%
}}}}
\put(19487,-19580){\makebox(0,0)[b]{\smash{{\SetFigFont{6}{7.2}{\rmdefault}{\mddefault}{\updefault}{\color[rgb]{0,0,0}$1$}%
}}}}
\put(18767,-19130){\makebox(0,0)[b]{\smash{{\SetFigFont{6}{7.2}{\rmdefault}{\mddefault}{\updefault}{\color[rgb]{0,0,0}$2$}%
}}}}
\put(19217,-19400){\makebox(0,0)[b]{\smash{{\SetFigFont{6}{7.2}{\rmdefault}{\mddefault}{\updefault}{\color[rgb]{0,0,0}$2$}%
}}}}
\put(19667,-19670){\makebox(0,0)[b]{\smash{{\SetFigFont{6}{7.2}{\rmdefault}{\mddefault}{\updefault}{\color[rgb]{0,0,0}$2$}%
}}}}
\put(18947,-19220){\makebox(0,0)[b]{\smash{{\SetFigFont{6}{7.2}{\rmdefault}{\mddefault}{\updefault}{\color[rgb]{0,0,0}$3$}%
}}}}
\put(19397,-19490){\makebox(0,0)[b]{\smash{{\SetFigFont{6}{7.2}{\rmdefault}{\mddefault}{\updefault}{\color[rgb]{0,0,0}$3$}%
}}}}
\put(18677,-19040){\makebox(0,0)[b]{\smash{{\SetFigFont{6}{7.2}{\rmdefault}{\mddefault}{\updefault}{\color[rgb]{0,0,0}$4$}%
}}}}
\put(19127,-19310){\makebox(0,0)[b]{\smash{{\SetFigFont{6}{7.2}{\rmdefault}{\mddefault}{\updefault}{\color[rgb]{0,0,0}$4$}%
}}}}
\put(19577,-19580){\makebox(0,0)[b]{\smash{{\SetFigFont{6}{7.2}{\rmdefault}{\mddefault}{\updefault}{\color[rgb]{0,0,0}$4$}%
}}}}
\put(18857,-19130){\makebox(0,0)[b]{\smash{{\SetFigFont{6}{7.2}{\rmdefault}{\mddefault}{\updefault}{\color[rgb]{0,0,0}$5$}%
}}}}
\put(19307,-19400){\makebox(0,0)[b]{\smash{{\SetFigFont{6}{7.2}{\rmdefault}{\mddefault}{\updefault}{\color[rgb]{0,0,0}$5$}%
}}}}
\end{picture}%